\input amstex
\documentstyle{amsppt}
\magnification1200
\pagewidth{6.5 true in}
\pageheight{9.25 true in}

\NoBlackBoxes
\topmatter
\title 
The number of imaginary quadratic fields with a given class number 
\endtitle
\author 
K. Soundararajan
\endauthor 
\address 
Department of Mathematics, Stanford University, Stanford, CA 94305, USA
\endaddress
\email 
ksound\@stanford.edu 
\endemail

\endtopmatter

\def\sumf{\sideset \and^\flat \to \sum} 
\document

\noindent Gauss asked for a list of imaginary quadratic fields with class number one.   
This problem inspired a great deal of work; some of the prominent milestones 
being the work of Heilbronn showing that there are only finitely many 
fields with a given class number, the work of Landau and Siegel providing good (but 
ineffective) lower bounds for the class number, the work of 
Heegner, Baker, and Stark showing that there are exactly nine fields with 
class number $1$, and the effective resolution of the class number problem 
due to Goldfeld, Gross and Zagier.   In this note we investigate the number, ${\Cal F}(h)$, of 
imaginary quadratic fields with class number $h$; thus, ${\Cal F}(1)=9$ is the 
celebrated Heegner-Baker-Stark result.   From Tatuzawa's refinement 
of the Landau-Siegel theorem one could 
compute ${\Cal F}(h)$ up to an error of $1$ relatively easily.  The Goldfeld-Gross-Zagier 
result permits, with great effort, the calculation of ${\Cal F}(h)$ for any given $h$, and Watkins [5]  
has accomplished this for all $h\le 100$.  What is the asymptotic behavior of ${\Cal F}(h)$ 
for large $h$?  This question is independent of the Landau-Siegel zero issue;
nevertheless it seems difficult to answer.  We establish here an asymptotic 
formula for the average value of ${\Cal F}(h)$, a modest non-trivial upper bound for 
${\Cal F}(h)$ (together with an application to a question of 
Rosen and Silverman on odd parts of class numbers), and we speculate on the nature of ${\Cal F}(h)$.

Throughout we let $-d$ denote a negative fundamental discriminant,  $\chi_{-d}$ 
will denote the associated primitive quadratic character $\pmod {|d|}$, and $h(-d)$ will 
denote the class number of ${\Bbb Q}(\sqrt{-d})$.  When $d>4$ recall that 
Dirichlet's class number formula gives 
$$
h(-d) = \sqrt{d} L(1,\chi_{-d})/\pi.
$$
Typically $L(1,\chi_{-d})$ has constant size; rarely does it fall 
outside the range $(1/10,10)$.  Therefore we would expect that class numbers 
below $H$ arise mainly from fields with discriminants of size about $H^2$, and 
the number of such fields should be asymptotically a constant times $H^2$.

\proclaim{Theorem 1}  As $H \to \infty$ we have 
$$
\sum_{h\le H} {\Cal F}(h) = \frac{3\zeta(2)}{\zeta(3)} H^2+O\Big(H^2 (\log H)^{-\frac 12+\epsilon}\Big).
$$
\endproclaim

\proclaim{Theorem 2} For large $H$ we have
$$
{\Cal F}(H)\ll H^2 \frac{(\log \log H)^4}{\log H}.
$$
\endproclaim

From Watkins [5] we know that there are $42272$ fields with class 
number below $100$; the main term of the asymptotic in  
Theorem 1 is approximately $41053$.  By modifying our argument 
one could improve the error term in the asymptotic 
formula of Theorem 1 to $O(H^2 (\log H)^{-1+\epsilon})$.   
Some new ideas seem needed to improve 
the power of $\log h$ appearing in Theorem 2.
 
\def\lam{\lambda}

We expect that ${\Cal F}(h)$ is of size about $h$ (the average size), although there is 
some variation.  More precisely, we conjecture that 
$$
\frac{h}{\log h} \ll {\Cal F}(h) \ll h \log h. \tag{C1} 
$$
Our heuristic reasoning is as follows.  Let $2^{\lambda}$ denote
the exact power of $2$ dividing $h$.  By genus theory, if the class 
number is $h$ then the fundamental discriminant $-d$ can have at
most $(\lam+1)$ prime factors if $-d \equiv 1 \pmod 4$, and 
$-d/4$ can have at most $(\lam+1)$ prime factors 
if $-d \equiv 0 \pmod 4$.  By the class number formula we also know that 
these discriminants are essentially of size $h^2$.  If $\ell \le \lam+1$ 
then there are $\asymp \frac{h^2}{\log h} \frac{(\log \log h)^{\ell-1}}{(\ell-1)!}$ 
fundamental discriminants of size $h^2$ with $-d$ (or $-d/4$) divisible by exactly 
$\ell$ primes.  For such discriminants the class number is of size about $h$, and 
constrained to be a multiple of $2^{\ell-1}$.  Thus we may think of the probability of 
the class number being exactly $h$ as roughly $2^{\ell-1}/h$.  In other words 
we expect that there are $\asymp 2^{\ell-1} h(\log \log h)^{\ell-1}/((\ell-1)! \log h)$ fields 
with $d$ (or $d/4$) composed of exactly $\ell$ prime factors, and with class 
number equal to $h$.   Summing over all $\ell \le \lam+1$ we 
arrive at 
$$
{\Cal F}(h)  \asymp \frac{h}{\log h} \sum_{\ell \le \lam+1}  \frac{2^{\ell-1} (\log \log h)^{\ell-1}}{(\ell-1)!}. 
\tag{C2}
$$

The unspecified constant in (C2) seems delicate, and would probably depend on 
arithmetical properties of $h$.  For example, the Cohen-Lenstra heuristics [1] predict 
that the probability of class numbers being divisible by $3$ is larger than $1/3$.  
So we would expect ${\Cal F}(h)$ to be larger when $3$ divides $h$, rather than 
when $3 \nmid h$.  Similar (smaller) biases would exist when $5$ divides $h$ etc.  
Inspecting Watkins' table (page 936 of [5]) we can already see the bias in favor 
of multiples of $3$.  

Conjecture (C2) does not lend itself to numerical testing.  To provide falsifiable 
conjectures, we may consider the ratio ${\Cal F}(h_1)/{\Cal F}(h_2)$  
for various choices of $h_1$ and $h_2$.  For example, if $h_1$ and $h_2$ are 
primes with $h_1/2 \le h_2 \le 2 h_1$ (say) then it seems safe to conjecture 
that 
$$
\frac{{\Cal F}(h_1)}{{\Cal F}(h_2)} \sim \frac{h_1}{h_2}. \tag{C3} 
$$
Also if $h$ is odd, and large then (C2) suggests that ${\Cal F}(h){\Cal F}(4h)/{\Cal F}(2h)^2$ 
should tend to $1/2$.  It would be interesting to assemble numerical data on these 
questions. 

This note was motivated by the recent paper of Rosen and Silverman [3] 
where they ask for information on $N(C;X)$ which counts the 
number of fundamental discriminants $-d$ with $1\le d\le X$ 
such that $h^{\text {odd}}(-d)$ (the odd part of the class number; in other words, the largest odd number 
dividing $h(-d)$) lies below a fixed number $C$.   Rosen and Silverman wished 
to know if $N(C;X) = o(X)$ for large $X$.   We show that such is indeed the case. 

\proclaim{Corollary 3}  For a fixed number $C$, and large $X$ we have 
$$
N(C;X) \ll X (\log \log X)^6/\log X.
$$
\endproclaim

\demo{Proof of Theorem 1}  We first show that one can restrict attention to discriminants 
$-d$ with $1\le d\le X=H^2(\log \log H)$.  First consider the range $X<d<H^2(\log H)^4$.  If $h(-d)<H$ then we must have $L(1,\chi_{-d}) \ll (\log \log H)^{-1/2}$, and by Theorem 4 of [2] 
there are at most\footnote{In fact Theorem 4 of [2] gives a much better bound, but the bound 
given above suffices for our purposes.}  $H^2/(\log H)$ values of $d<H^2(\log H)^4$ with such a small 
value of $L(1,\chi_{-d})$.   
If $d>H^2(\log H)^4$ then $h(-d)$ can be below $H$ only 
when $L(1,\chi_{-d}) \ll H/\sqrt{d}\ \  (\le 1/(\log H)^2)$.  
Tatuzawa's theorem (see [4])  
shows that there is at most one such discriminant $-d$ with $d>H^2 (\log H)^4$.  
Therefore 
\def\sumf{\sideset \and^{\flat} \to \sum} 
$$
\sum\Sb h\le H\endSb {\Cal F}(h)  = \sumf\Sb d\le X\\ h(-d)\le H\endSb 1 + O\Big(\frac{H^2}{\log H}\Big), 
\tag{1}
$$
where the $\flat$ indicates that the sum is over fundamental discriminants $-d$.  

Observe that for any $c>0$, 
$$
\frac{1}{2\pi i} \int_{c-i\infty}^{c+i\infty} 
\frac{x^s}{s} \Big(\frac{(1+\delta)^{s+1}-1}{\delta(s+1)} \Big) 
ds 
= \cases
1 &\text{if } x \ge 1 \\ 
(1+\delta-1/x)/\delta &\text{if  } (1+\delta)^{-1} \le x\le 1\\ 
0&\text{if } x\le (1+\delta)^{-1}.\\ 
\endcases
$$
Here $\delta >0$ is a parameter which we shall choose later.  By the 
class number formula and (1) we get that 
$$
\align
\sum_{h\le H} {\Cal F}(h)& \le  
\frac{1}{2\pi i} \int_{c-i\infty}^{c+i\infty}  
\sumf\Sb d\le X \endSb \Big(\frac{\pi }{\sqrt{d} L(1,\chi_{-d})}\Big)^{s} \frac{H^s}{s} \Big( 
\frac{(1+\delta)^{s+1}-1}{\delta (s+1)} \Big) ds + O\Big(\frac{H^2}{\log H}\Big)\\
&\le 
\sum_{h\le H(1+\delta)}  {\Cal F}(h).
\tag{2}  \\
 \endalign
 $$
 We now focus on evaluating the integral in (2) which leads naturally to Theorem 1.  
 
 We shall take $c=1/\log H$, and $\delta =({\log H})^{-\frac 12}$.   Set $S=\log X/(10^4 
 (\log \log X)^2)$.  The region 
 $|s| > S$ contributes to the integral in (2) an amount 
 $$
 \ll \frac{X}{\delta } \int_{|s|> S} \frac{1}{|s(s+1)|} |ds| 
 \ll H^2 (\log H)^{-\frac 12+\epsilon}. \tag{3}
 $$
 In the region $|s| \le  S$ we shall use Theorem 2 of 
 [2] in order to evaluate the sum over $d$.  That result evaluates such sums in terms of 
 a probabilistic model for $L(1,\chi_{-d})$.

 For primes $p$ let ${\Bbb X}(p)$ denote 
 independent random variables taking the value $1$ with probability 
 $p/(2(p+1))$, $0$ with probability $1/(p+1)$, and $-1$ with 
 probability $p/(2(p+1))$.  Let $L(1,{\Bbb X}) = \prod_{p} (1-{\Bbb X}(p)/p)^{-1}$.  
 This product converges almost surely, and the main results of [2] 
 compare the distribution of $L(1,\chi_{-d})$ with the distribution of 
 such random Euler products.   With two caveats that we clarify below, 
 Theorem 2 of [2] gives that 
 for $|z|\le \log x/(500 (\log \log x)^2)$ and Re$(z)>-1$ 
 $$
 \sumf\Sb d\le x\endSb L(1,\chi_{-d})^z = \frac{3}{\pi^2}x {\Bbb E}(L(1,{\Bbb X})^z) 
 + O\Big(x \exp \Big(-\frac{\log x}{5\log \log x}\Big)\Big), \tag{4}
 $$
 where ${\Bbb E}$ stands for expectation. 
 The first caveat is that Theorem 2 of [2] considers both positive and negative 
 fundamental discriminants, but the arguments given there permit us to 
 restrict to negative fundamental discriminants as above.  Secondly, there we omitted 
 a small number ($\ll \log x$) of exceptional Landau-Siegel discriminants.  Since 
 $L(1,\chi_{-d}) \gg 1/\sqrt{x}$ and $\text{Re}(z) \ge -1$ the contribution of 
 these exceptional discriminants to our sum is $\ll \sqrt{x} \log x$, and so 
 (4) holds.  Using (4) and partial summation we obtain that for 
 $|s|\le S$ and Re$(s)=1/\log H$ we have 
 $$
 \sumf_{d\le X} (\sqrt{d}L(1,\chi_{-d}))^{-s} = 
 \frac{3}{\pi^2} {\Bbb E}(L(1,{\Bbb X})^{-s}) \int_1^X x^{-s/2}dx + O\Big(X\exp \Big(-\frac{\log X}{5\log \log X}
 \Big)\Big). \tag{5}
 $$
 
 From (3) and (5) we see that the integral in (2) is, with an error $O(H^2 (\log H)^{-\frac 12+\epsilon})$, 
 $$
 \frac{1}{2\pi i} \int_{|s|\le S} 
 \frac{3}{\pi^2} {\Bbb E}(L(1,{\Bbb X})^{-s}) \Big(\int_1^X x^{-s/2} dx \Big) \frac{(\pi H)^s}{s} 
 \Big(\frac{(1+\delta)^{s+1}-1}{\delta(s+1)} \Big) ds. \tag{6}
  $$
 For $1\le x\le X$ we may see that 
 $$
 \align
 \frac{1}{2\pi i}  \int_{|s| \le S} 
 \Big( \frac{\pi H}{\sqrt{x}L(1,{\Bbb X})}\Big)^s \frac 1s &\Big(\frac{(1+\delta)^{s+1}-1}{\delta(s+1)} 
 \Big) ds
= O\Big(\frac{L(1,{\Bbb X})^{-c}}{(\log H)^{\frac 12-\epsilon}}\Big) 
\\
&+ \cases 
1 &\text{if } \sqrt{x} L(1,{\Bbb X}) <\pi H \\ 
\in [0,1] &\text{if }\pi H\le \sqrt{x} L(1,{\Bbb X})\le \pi H(1+\delta)\\ 
0 &\text{if } \pi H(1+\delta) < \sqrt{x} L(1,{\Bbb X}). 
\endcases\\
\endalign
$$
Integrating this over $x$ from $1$ to $X$ we get
$$
O\Big( \frac{H^2}{(\log H)^{\frac 12-\epsilon}} (1+L(1,{\Bbb X})^{-c})\Big) 
+ \min\Big( \frac{\pi^2 H^2}{L(1,{\Bbb X})^2}, X\Big). 
$$
Therefore the quantity in (6) equals 
$$
{\Bbb E} \Big( \min\Big( \frac{3H^2}{L(1,{\Bbb X})^2}, \frac{3X}{\pi^2}\Big) \Big) 
+ O\Big( \frac{H^2}{(\log H)^{\frac 12-\epsilon}}\Big), \tag{7}
$$
and this is also our integral in (2).

Proposition 1 of [2] reveals that the probability that $L(1,{\Bbb X})$ is 
less than $\pi^2/(6e^{\gamma}\tau)$ is $\exp(-e^{\tau-C_1}/\tau +O(e^{\tau}/\tau^2))$ for 
some absolute constant $C_1$.  
Hence we may see that 
$$
{\Bbb E}\Big( \min\Big( \frac{3H^2}{L(1,{\Bbb X})^2}, \frac{3X}{\pi^2}\Big) \Big) 
={\Bbb E}\Big(\frac{3H^2}{L(1,{\Bbb X})^2}\Big) + O\Big( \frac{H}{\log H}\Big). 
$$
Finally, by independence of the random variables ${\Bbb X}(p)$ we 
have 
$$
\align
{\Bbb E}(L(1,{\Bbb X})^{-2}) &= \prod_{p} {\Bbb E}\Big( \Big( 1-\frac{{\Bbb X}(p)}{p}\Big)^2\Big) 
\\
&= \prod_{p} \Big(\frac{p}{2(p+1)} \Big(1-\frac 1p\Big)^2 + \frac{1}{(p+1)} + \frac{p}{2(p+1)} \Big(1+\frac 1p\Big)^2 \Big)\\
&=\prod_p \Big(1-\frac{1}{p^3}\Big) \Big(1-\frac{1}{p^2}\Big)^{-1} = 
\frac{\zeta(2)}{\zeta(3)}.\\
\endalign
$$
Using these observations in (7), we conclude that the integral in (2) is 
$$
\frac{3\zeta(2)}{\zeta(3)} H^2 + O \Big(\frac{H^2}{(\log H)^{\frac 12-\epsilon}}\Big).
$$
This establishes Theorem 1.

\enddemo

\demo{Proof of Theorem 2} As before
set $X=H^2 \log \log H$, and $S= (\log X)/(10^4(\log \log X)^2)$.  As in (1) we see that 
$$
{\Cal F}(H) = \sumf\Sb d\le X \\ h(-d)=H\endSb 1 +O \Big(\frac{H^2}{\log H}\Big).
$$
Since 
$$
\frac{1}{S} \int_{-S}^{S} \Big(1-\frac{|x|}{S}\Big) e^{2\pi i x\xi} dx 
\qquad
\cases
=1 &\text{if } \xi = 0,\\ 
\ge 0 &\text{always},\\
\endcases 
$$
we deduce, by the class number formula, that 
$$
{\Cal F}(H) \le O\Big (\frac{H^2}{\log H}\Big) + \frac{1}{S} \int_{-S}^{S} 
\Big(1-\frac{|x|}{S}\Big) \sumf_{d\le X} \Big(\frac{\pi H}{\sqrt{d}L(1,\chi_{-d})}\Big)^{ix} dx.\tag{8}
$$
As in (5) we have that 
$$
\sumf_{d\le X} (\sqrt{d} L(1,\chi_{-d})^{-ix} 
= \frac{3}{\pi^2} {\Bbb E}(L(1,{\Bbb X})^{-ix}) \int_1^X y^{-ix/2}dy + O\Big(\frac{H}{(\log H)^2}\Big)
\ll \frac{X}{1+|x|} + \frac{H}{(\log H)^2}.
$$
Inserting this in (8) we obtain that 
$$
{\Cal F}(H) \ll \frac{H^2}{\log H} + X \frac{\log S}{S} \ll H^2 \frac{(\log \log H)^4}{\log H}.
$$

\enddemo

\demo{Proof of Corollary 3}   From Theorem 4 of [2] (with $\tau$ there being $\log \log X$) 
we have that the number of fundamental discriminants $-d$ with $1\le d\le X$ and $h(-d) > \sqrt{X} \log \log X$ is at most $X\exp(-c\log X/\log \log X)$ for some positive constant $c$.  
Therefore 
$$
N(C,X) \le X\exp\Big(-c\frac{\log X}{\log \log X}\Big) + \sum\Sb 2^k \ell \le \sqrt{X}\log \log X\\ 
\ell \text{ odd}\\ \ell \le C\endSb {\Cal F}(2^k\ell).
$$
The Corollary now follows from Theorem 2.
\enddemo

\enddocument